\documentclass[12pt]{amsart}

\usepackage[T1]{fontenc}
\usepackage[utf8]{inputenc}
\usepackage{amsaddr}

\usepackage{amsmath,amssymb,amsthm,mathrsfs}
\usepackage{url}
\usepackage{enumitem}
\usepackage[colorlinks=true,linkcolor=blue,citecolor=blue]{hyperref}
\numberwithin{equation}{section}
\newtheorem{thm}{Theorem}[section]
\newtheorem{prop}[thm]{Proposition}
\newtheorem{lem}[thm]{Lemma}
\newtheorem{cor}[thm]{Corollary}
\newtheorem*{mainthm}{Main Theorem}
\theoremstyle{definition}
\newtheorem{example}[thm]{Example}
\theoremstyle{remark}

\DeclareMathOperator{\GKdim}{GKdim}
\newcommand{\GK}{\GKdim}
\newcommand{\de}{\delta}
\newcommand{\ad}{\operatorname{ad}}
\DeclareMathOperator{\End}{End}

\title[GK dimension of constants of algebraic $q$-skew derivations]{Gelfand--Kirillov dimension of  constants of algebraic
$q$-skew derivations}
\author{Piotr Grzeszczuk}

\address{
Faculty of Computer Science\\
Bialystok University of Technology\\
Wiejska 45A\\
15-351 Białystok\\
Poland
}

\email{p.grzeszczuk@pb.edu.pl}

\thanks{This work was supported by Grant WZ/WI-IIT/2/2025 from the Bialystok University of Technology and funded by the Ministry of Science and Higher Education of Poland through resources allocated for research.}


\begin{document}

\begin{abstract}
Let $\de$ be an algebraic $q$-skew $\sigma$-derivation of an algebra $R$ over
a field of arbitrary characteristic, where $\sigma$ is an algebraic
automorphism of $R$, and let $R^{\de}$ be the algebra of constants of $\de$.
We show that $\GK(R) = \GK(R^{\de})$ whenever $R$ is semiprime and $R^{\de}$
is a PI algebra containing no infinite direct sum of nonzero ideals. The polynomial identity
cannot be dropped. For every $r \ge 1$ we construct a primitive algebra $R$
with an inner derivation $d$ such that $d^{3} = 0$, $\GK(R) \ge 2r$ and
$\GK(R^{d}) = 0$, while $R^{d}$ contains no direct sum of two nonzero ideals.
\end{abstract}

\maketitle

\section{Introduction}\label{sec:intro}

Throughout, $R$ is an associative algebra over a field $F$ and $\GK$ denotes
the Gelfand--Kirillov dimension over $F$.  Let $\sigma$ be an $F$-linear
automorphism of $R$.  Recall that an $F$-linear map $\de\colon R \to R$ is a
\emph{$\sigma$-derivation} if
$\de(rs) = \de(r)s + \sigma(r)\de(s)$ for all $r,s \in R$.  We say that $\de$
is \emph{$q$-skew} if $\de\sigma = q\,\sigma\de$ for some nonzero scalar $q \in F^\times$.
The \emph{algebra of constants} of $\de$ is the subalgebra
$R^{\de} = \ker\de$.  Derivations correspond to $\sigma = \mathrm{id}$, and we
denote them by $d$.  A linear map is called \emph{algebraic} if it satisfies
a nonzero polynomial over $F$.

A natural question is whether $R$ and $R^{d}$ have the same growth.  In
general they do not, even for locally nilpotent derivations.  Lorenz
\cite[Example]{L} showed that for every $r \ge 1$ there is a locally
finite-dimensional algebra $R$ with a derivation $d$ such that the
differential operator ring $T = R[x;d]$ is affine and $\GK(T) = r$, while
$\GK(R) = 0$.  He also gave an example with $\GK(T) = \infty$.  The algebra
$T$ has a locally nilpotent derivation $\partial = d/dx$, given by
$\partial(R) = 0$ and $\partial(x) = 1$.  In characteristic $0$ its algebra
of constants is $R$.  Hence $\GK(T) - \GK(T^{\partial})$ can be arbitrarily
large, or even infinite.  By contrast, $\GK(T) = \GK(R) + 1$ whenever every
finite-dimensional subspace of $R$ is contained in a $d$-stable finitely
generated subalgebra of $R$ \cite[Proposition 3.5]{KL}.

Under additional assumptions the difference is at most $1$.  In
\cite[Theorem 10]{BG3} it is proved that if $d$ is a locally nilpotent
derivation of a finitely generated semiprime PI algebra $A$ of
characteristic $0$ and finite GK dimension, then $\GK(A) = \GK(A^{d})$
when $d$ is nilpotent on an ideal of $A$ with the same GK dimension as
$A$, and $\GK(A) = \GK(A^{d}) + 1$ otherwise.  In particular, the equality holds for nilpotent derivations.  This suggests that, under suitable hypotheses, the equality $\GK(R) = \GK(R^{d})$ holds for algebraic derivations.  We prove this in the more general setting of skew derivations.  Our main result is the following.

\begin{mainthm}
Let $R$ be a semiprime algebra over a field of arbitrary characteristic, let
$\sigma$ be an automorphism of $R$ and let $\de$ be a $q$-skew
$\sigma$-derivation of $R$, with $\de$ and $\sigma$ algebraic.  If $R^{\de}$
satisfies a polynomial identity and contains no infinite direct sum of nonzero ideals,
then $\GK(R) = \GK(R^{\de})$.
\end{mainthm}

The proof first lifts both hypotheses from $R^{\de}$ to $R$: finite prime
dimension through the generalized $0$-eigenspace $R_0$ of $\de$, by Lemma
\ref{lem:primedim} and \cite[Proposition 4.4]{BG00}, and the polynomial
identity by Theorem \ref{thm:PIsemiprime}, which extends
\cite[Theorem 7.1]{BG00} beyond prime algebras and the remark after
\cite[Theorem 6.3]{BG00} beyond roots of unity $q$.  We
then pass to the centre.  By Lemma \ref{lem:growthcentre},
$\GK(R) = \GK(Z(R))$, and $\sigma$ has finite order on $Z(R)$, so it suffices
to show that $Z(R)^{\sigma}$, or in characteristic $p$ the set of $p$-th
powers of its elements, lies in $R^{\de}$ (Section \ref{sec:qskew}).

Section \ref{sec:examples} shows that the hypotheses cannot be substantially
weakened.  Without  finite prime dimension of
$R^{d}$ and semiprimeness of $R$, the difference $\GK(R) - \GK(R^{d})$ can be arbitrarily large, or
infinite, with $d^{2} = 0$ and $R^{d}$ commutative of GK dimension $0$
(Example \ref{ex:biggap}).  Semiprimeness of $R$ together with finite
prime dimension of $R^{d}$ does not suffice either, even for primitive $R$:
for every $r \ge 1$ there is a primitive
algebra $R$ with an inner derivation $d$, $d^{\,3} = 0$, such that
$\GK(R) \ge 2r$ and $\GK(R^{d}) = 0$, while $R^{d}$ contains no direct sum
of two nonzero ideals.  Here $R^{d}$ is not PI, so the polynomial identity
cannot be dropped (Example \ref{ex:primitive}).

\section{Preliminaries}\label{sec:prelim}

We use the basic properties of Gelfand--Kirillov dimension from \cite{KL}:
monotonicity on subalgebras, including non-unital ones, and the
computation of $\GK$ as the supremum over affine subalgebras.

\begin{prop}[{\cite[Proposition 5.5]{KL}, \cite[8.2.9]{MR}}]\label{prop:bimod}
Let $B \subseteq R$ be a subalgebra such that $R$ is finitely generated as a
left $B$-module.  Then $\GK(R) = \GK(B)$.
\end{prop}

\begin{cor}\label{cor:integral}
Let $A$ be a commutative algebra and $B \subseteq A$ a subalgebra, and suppose
that $A$ is spanned by elements each of which satisfies a monic equation of
degree at most $t$ over $B$.  Then $\GK(A) = \GK(B)$.
\end{cor}

\begin{proof}
Only $\le$ needs proof.  A finite-dimensional $V \subseteq A$ lies in the span
of finitely many of the given elements, say $v_1,\dots,v_k$.  Let
$B_0 \subseteq B$ be the affine subalgebra generated by the coefficients of
monic equations of degree at most $t$ for the $v_i$.  As $A$ is commutative,
reducing exponents modulo those equations shows that the algebra generated by
$v_1,\dots,v_k$ and $B_0$ is generated as a $B_0$-module by the $t^{\,k}$
monomials $v_1^{\,b_1}\cdots v_k^{\,b_k}$ with all $b_i < t$.  Proposition \ref{prop:bimod} gives
$\GK(F[V \cup B_0]) = \GK(B_0) \le \GK(B)$, and taking the supremum over $V$
yields $\GK(A) \le \GK(B)$.
\end{proof}

A ring is \emph{semiprime} if it has no nonzero nilpotent ideals. An algebra has \emph{finite prime dimension} if it contains no infinite direct
sum of nonzero two-sided ideals.  For a semiprime algebra $R$ the
following conditions are equivalent to finite prime dimension:
\begin{itemize}
\item $R$ has finite Goldie rank as an $R \otimes_F R^{op}$-module;
\item $R$ contains ideals $I_1,\dots,I_k$ whose sum is direct and essential,
each $I_i$ being a prime ring;
\item $R$ has finitely many minimal prime ideals $P_1,\dots,P_k$ with
$\bigcap_i P_i = 0$.
\end{itemize}
The equivalence follows from \cite[Theorem 2.2.15]{MR}, 
and we use all three
forms, related by $P_i = \operatorname{ann}(I_i)$ for
$I_i = \bigcap_{j \neq i}P_j$.

Finite prime dimension passes to ideals of a semiprime $R$.  Indeed, if $A$
is an ideal of $R$ and $B \neq 0$ an ideal of the ring $A$, then $ABA
\subseteq B$ is an ideal of $R$, and $ABA \neq 0$: otherwise
$(BA)^{2} = B(ABA) = 0$, so the right ideal $BA$ of $R$ is zero.  Thus $B$
lies in the left annihilator $\{x \in A : xA = 0\}$ of $A$ in $A$, which is
an ideal of $R$ of square zero, hence zero, as $R$ is semiprime.  Consequently, an infinite direct sum of nonzero ideals of $A$ gives one in $R$.  This argument does
not use a unit, so it applies to ideals $A$ regarded as rings without unit.

\begin{lem}\label{lem:primedim}
Let $A$ be an algebra, $\sigma$ an automorphism and $\de$ a $q$-skew
$\sigma$-derivation of $A$ with $\de^{\,n} = 0$, $n \ge 1$.  If $A^{\de}$
contains no infinite direct sum of nonzero ideals, then $A$ contains no
infinite direct sum of nonzero ideals, that is, $A$ has finite prime
dimension.
\end{lem}

\begin{proof}
Put $\mathcal{F}_{k} = \ker\de^{\,k+1}$, so that $\mathcal{F}_{0} = A^{\de}$
and $\mathcal{F}_{n-1} = A$.  As $q \neq 0$, from
$\de^{\,k}\sigma = q^{k}\sigma\de^{\,k}$ and
$\de^{\,k}\sigma^{-1} = q^{-k}\sigma^{-1}\de^{\,k}$ we get
$\sigma(\mathcal{F}_k) = \mathcal{F}_k$ for every $k$, and in particular
$\sigma(A^{\de}) = A^{\de}$.  For $b,c \in A^{\de}$ and $a \in A$, we have $\de(bac) = \de(b)ac + \sigma(b)\de(ac) = \sigma(b)\de(a)c$.  As $\sigma(b) \in A^{\de}$, induction gives
$\de^{\,j}(bac) = \sigma^{j}(b)\,\de^{\,j}(a)\,c$ for all $j$, so each
$\mathcal{F}_k$ is an $A^{\de}$-sub-bimodule of $A$.
For an $A^{\de}$-bimodule $M$ let ${}^{\sigma}\!M$ be $M$ with the left action
twisted by $\sigma$, that is $b \cdot m = \sigma(b)m$.  As $\sigma$ is an
automorphism of $A^{\de}$, the bimodules $M$ and ${}^{\sigma}\!M$ have the same
subbimodules.  By the formula above, $\de$
restricts to a homomorphism of $A^{\de}$-bimodules
$\mathcal{F}_k \to {}^{\sigma}\mathcal{F}_{k-1}$ with kernel $A^{\de}$, so
that there are left exact sequences
$$
0 \;\longrightarrow\; A^{\de} \;\longrightarrow\; \mathcal{F}_{k}
\;\stackrel{\de}{\longrightarrow}\; {}^{\sigma}\mathcal{F}_{k-1}
\qquad (1 \le k \le n-1) .
$$
We regard them as sequences of modules over
$S = A^{\de} \otimes_F (A^{\de})^{\mathrm{op}}$. The $S$-submodules of
$A^{\de}$ are precisely the two-sided ideals of $A^{\de}$. By hypothesis,
$A^{\de}$ has finite uniform dimension $s$ as an $S$-module. Since uniform
dimension is subadditive on extensions, monotone on submodules, and invariant
under twisting by $\sigma$, we obtain
$$
\operatorname{udim}_S(\mathcal{F}_{k}) \;\le\; \operatorname{udim}_S(A^{\de})
+ \operatorname{udim}_S\bigl(\de(\mathcal{F}_{k})\bigr)
\;\le\; s + \operatorname{udim}_S(\mathcal{F}_{k-1}).
$$
By induction, $\operatorname{udim}_S(A) \le ns$. Every ideal of $A$ is an
$A^{\de}$-sub-bimodule of $A$, and hence an $S$-submodule. Consequently, $A$
cannot contain an infinite direct sum of nonzero ideals.
\end{proof}

\begin{lem}\label{lem:growthcentre}
Let $A$ be a semiprime PI algebra with finitely many minimal prime ideals.
Then $\GK(A) = \GK\bigl(Z(A)\bigr)$.
\end{lem}

\begin{proof}
Only the inequality $\le$ requires proof. Let $P_1,\dots,P_k$ be the minimal
prime ideals of $A$, and for each $i$ let $\pi_i \colon A \to A_i = A/P_i$ be
the canonical projection and $J_i = \bigcap_{j \neq i}P_j$. Since $A$ embeds
in $\prod_i A_i$, we have $\GK(A) = \max_i \GK(A_i)$. Moreover, $\GK(A_i) = \GK(Z_i)$, where $Z_i = Z(A_i)$.  Indeed, since $A_i$
is prime PI, by Rowen's theorem \cite{Ro} and Posner's theorem \cite{Po} the
localization $A_iZ_i^{-1}$ at the nonzero central elements is a simple
algebra, finite dimensional over its centre, the field of fractions
$Z_iZ_i^{-1}$ of $Z_i$.  Central localization does not change GK dimension
\cite[Chapter 4]{KL}, so Proposition \ref{prop:bimod} gives
$\GK(A_i) \le \GK(A_iZ_i^{-1}) = \GK(Z_iZ_i^{-1}) = \GK(Z_i) \le \GK(A_i)$.

Fix $i$. Then $J_i \neq 0$, by minimality of $P_i$, and $J_i \cap P_i = 0$, so $\pi_i$ maps $J_i$
isomorphically onto a nonzero ideal of $A_i$. If $x \in J_i$ and $\pi_i(x)$
is central in $A_i$, then every commutator $[x,a]$ with $a \in A$ lies in
$J_i \cap P_i = 0$, and hence $x \in Z(A)$. It follows that $\pi_i$ maps
$J_i \cap Z(A)$ isomorphically onto $I = \pi_i(J_i) \cap Z(A_i)$, which is an
ideal of the domain $D = Z(A_i)$ and is nonzero by Rowen's theorem
\cite{Ro}. So $\GK(Z(A)) \ge \GK(I)$, and $\GK(I) \ge \GK(D)$: for
$0 \neq c \in I$ and an affine $B = F[z_1,\dots,z_s] \subseteq D$, the space
$V = \operatorname{span}\{c, cz_1,\dots,cz_s\} \subseteq I$ satisfies
$V^{n+1} \supseteq c^{\,n+1}B_n$, where $B_n$ is the span of the monomials
of degree at most $n$, and multiplication by $c^{\,n+1}$ is injective.
Hence $\dim V^{n+1} \ge \dim B_n$, so $\GK(I) \ge \GK(B)$, and therefore
$\GK(A) = \max_i \GK(Z_i) \le \GK(Z(A))$.
\end{proof}

\section{\texorpdfstring{Algebraic $q$-skew $\sigma$-derivations}{Algebraic q-skew sigma-derivations}}\label{sec:qskew}

In this section we prove the Main Theorem (as Theorem \ref{thm:qskew}). Let $\de$ be a
$q$-skew $\sigma$-derivation of $R$ and write $R_0 = \{r \in R : \de^{\,k}(r) = 0 \text{ for some } k\}$ for the generalized $0$-eigenspace of $\de$.  It is a subalgebra.  Indeed, the
$q$-Leibniz rule writes $\de^{\,n}(xy)$ as a combination of the products
$\sigma^{k}\bigl(\de^{\,n-k}(x)\bigr)\de^{\,k}(y)$, $0 \le k \le n$, and
$\sigma$ maps $\ker\de^{\,j}$ into itself, as $\de\sigma = q\,\sigma\de$.  So if
$\de^{\,a}(x) = 0 = \de^{\,b}(y)$, then $\de^{\,a+b}(xy) = 0$.  The map
$\de$ acts locally nilpotently on $R_0$, and nilpotently because $\de$ is
algebraic, and $R_0^{\,\de} = R^{\de}$.

We collect the results of \cite{BG00} which we use.

\begin{thm}[\cite{BG00}]\label{thm:bg00}
Let $\de$ be a $q$-skew $\sigma$-derivation of a semiprime algebra $R$, with
$\de$ and $\sigma$ algebraic.
\begin{enumerate}
\item $R$ has finite prime dimension if and only if $R_0$ has, and if $R_0$
satisfies a polynomial identity, then so does $R$ \cite[Proposition
4.4]{BG00}.
\item If $R$ is prime and $R^{\de}$ satisfies a polynomial identity, then so
does $R$ \cite[Theorem 7.1]{BG00}.
\item If $R$ has finite prime dimension, $q$ is a root of unity and $R^{\de}$
satisfies a polynomial identity, then so does $R$ \cite[remark after
Theorem 6.3]{BG00}.  For $\sigma = \mathrm{id}$ this is
\cite[Proposition 5.4]{BG00}.
\end{enumerate}
\end{thm}

We need (3) without the restriction on $q$.

\begin{thm}\label{thm:PIsemiprime}
Let $\de$ be a $q$-skew $\sigma$-derivation of a semiprime algebra $R$ of
finite prime dimension, with $\de$ and $\sigma$ algebraic.  If $R^{\de}$
satisfies a polynomial identity, then so does $R$.
\end{thm}

\begin{proof}
In this proof, algebras are not assumed to be unital. Following
\cite{BG00}, we apply the results of that paper to nonzero ideals of
semiprime algebras. Recall that in a semiprime algebra every nilpotent
one-sided ideal is zero, and the left and right annihilators of any ideal $K$
coincide. We denote this common annihilator by $\operatorname{ann}(K)$. For an
automorphism $\tau$ of an algebra $B$ and $\lambda \in F$, set
$B_{\lambda} = \bigcup_{n \ge 1}\ker(\tau-\lambda)^{n}$. It is well known
that the sum of the subspaces $B_{\lambda}$ is direct and that
$B_{\lambda}B_{\mu} \subseteq B_{\lambda\mu}$ for all $\lambda, \mu \in F$.

\emph{Step 1.}  If $q$ is a root of unity, the assertion follows from
Theorem~\ref{thm:bg00}(3). Assume therefore that $q$ is not a
root of unity. After extending scalars to an algebraic closure $\overline{F}$
of $F$, the identity $(\de - q\lambda)^{n}\sigma = q^{n}\sigma(\de-\lambda)^{n}$ shows that $\sigma$ maps the generalized $\lambda$-eigenspace of $\de$ onto
the generalized $q\lambda$-eigenspace. Hence, if $\lambda \neq 0$ is an
eigenvalue of $\de$, then so is $q^{j}\lambda$ for every $j \ge 0$. Since
these elements are pairwise distinct and $\de$ has only finitely many
eigenvalues, $0$ is the only eigenvalue of $\de$, and therefore $\de$ is
nilpotent.

\emph{Step 2: reduction to the $\sigma$-prime case.}  The automorphism
$\sigma$ permutes the finitely many minimal primes of $R$.  For a
$\sigma$-orbit $O_j$ put $I_j = \bigcap_{P \notin O_j}P$.  The sum
$I_1 \oplus \dots \oplus I_h$ is an essential ideal, $\sigma(I_j) = I_j$,
$I_j = \operatorname{ann}\bigl(\bigoplus_{l\neq j}I_l\bigr)$ (if
$xI_l = 0$, then $xRI_l = 0$ and $I_l \not\subseteq P$ give $x \in P$ for
$P \in O_l$), and $I_j$ has finite prime dimension
(Section~\ref{sec:prelim}).

$I_j$ is $\sigma$-prime.  Indeed, let $A,B$ be $\sigma$-stable ideals of the ring
$I_j$ with $AB = 0$.  By Section~\ref{sec:prelim}, $A' = I_jAI_j \subseteq A$
and $B' = I_jBI_j \subseteq B$ are $\sigma$-stable ideals of $R$, nonzero if
$A$, resp.\ $B$, is, and $A'B' = 0$.  The set $X$ of $P \in O_j$ with
$A' \subseteq P$ is $\sigma$-stable, so $X = O_j$ or $X = \emptyset$, and in
the latter case $B' \subseteq P$ for all $P \in O_j$.  As $I_j \subseteq P$
for $P \notin O_j$, we get $A' = 0$ or $B' = 0$.

The annihilator of a $\sigma$-stable ideal $K$ is $\de$-stable.  Let
$x \in \operatorname{ann}(K)$ and $k,k' \in K$.  Applying $\de$ to $kx = 0$
gives $\sigma(k)\de(x)k' = -\de(k)xk' = 0$, so $K\de(x)K = 0$.  Then the
right ideal $\de(x)K$ and the left ideal $K\de(x)$ have square zero, hence
vanish.  So $\de(I_j) \subseteq I_j$, and $\de|_{I_j}$ is a nilpotent
$q$-skew $\sigma$-derivation with $I_j^{\,\de} \subseteq R^{\de}$.  If each
$I_j$ is PI, so is $\bigoplus_j I_j$, and so is $R$: a semiprime algebra
satisfies the multilinear identities of any essential PI ideal $J$, since
$R \subseteq Q_s(R) = Q_s(J)$ \cite[Theorem 2]{Ma}.  Therefore we may assume that $R$ is $\sigma$-prime.

\emph{Step 3: the element $a$.}  The minimal primes $P_1,\dots,P_k$ form one
$\sigma$-orbit, since for a proper nonempty $\sigma$-stable set of them the
intersections of the primes inside and outside it would be nonzero
$\sigma$-stable ideals with zero product.  If $k = 1$, then $R$ is prime and
we apply Theorem \ref{thm:bg00}(2).  Let $k \ge 2$.  Then $R$ is $\sigma$-prime, semiprime, with $k$ minimal
primes, $\de$ is nilpotent by Step~1 and $q$ is not a root of unity, so
\cite[Lemma 5.1]{BG00} applies (as in \cite{BG00}, to an algebra possibly
without unit): there are a nilpotent $a \in Q_s(R)$, say $a^{t} = 0$, and a
nonzero $\sigma$-stable ideal $I$ with $aI + Ia \subseteq I$ and
$\de(r) = ar - \sigma(r)a$ for $r \in R$.  As $\operatorname{ann}(I)$ is
$\sigma$-stable and $\operatorname{ann}(I)I = 0$, we have
$\operatorname{ann}(I) = 0$, so $I$ is essential.

Expanding $\de\sigma(r) = q\,\sigma\de(r)$ shows that $u = a - q\sigma(a)$
satisfies $ur = \sigma(r)u$ for $r \in R$, and then $uy = \sigma(y)u$ for
$y \in Q_s(R)$: if $Dy \subseteq R$ for an essential ideal $D$, then
$$\sigma(d)(uy - \sigma(y)u) = u(dy) - \sigma(dy)u = 0$$
 for $d \in D$.  Let
$e_i \in Q_s(R)$ be the projection of the essential ideal
$\bigoplus_i \bigcap_{j\neq i}P_j$ onto its $i$-th summand.  The $e_i$ are
orthogonal central idempotents with $\sum_i e_i = 1$, and
$\sigma(e_i) = e_{\pi(i)}$ with $\pi(i) \neq i$.  Hence
$e_iu = ue_i = e_{\pi(i)}u$ and $e_iu = e_ie_{\pi(i)}u = 0$, so $u = 0$, that
is, $\sigma(a) = q^{-1}a$.

The ideal $I$ is semiprime of finite prime dimension
(Section~\ref{sec:prelim}), $\sigma$- and $\de$-stable, as
$\de(I) \subseteq aI + Ia$, and $I^{\,\de} \subseteq R^{\de}$.  By Step 2,
$R$ is PI if $I$ is, and $Q_s(I) = Q_s(R)$.  Replacing $R$ by $I$ we may
assume, dropping $\sigma$-primeness,
$$
aR + Ra \subseteq R,\qquad \sigma(a) = q^{-1}a,\qquad a^{t} = 0,\qquad
\de(r) = ar - \sigma(r)a .
$$

\emph{Step 4: the subalgebra $T$.} For $\lambda \in F$, let $R_{\lambda}$
denote the generalized $\lambda$-eigenspace of $\sigma$ on $R$. For
$i \in \mathbb{Z}$ put $T_i = R_{q^{i}}$, and set $T = \bigoplus_{i\in\mathbb{Z}} T_i$.  Since $q$ is not a root of unity, the elements $q^{i}$ are pairwise
distinct, and since $R_{\lambda}R_{\mu} \subseteq R_{\lambda\mu}$, we have
$T_iT_j \subseteq T_{i+j}$. Thus $T$ is a $\sigma$-stable
$\mathbb{Z}$-graded subalgebra of $R$ with $T_0 = R_1$. Since
$\sigma(a) = q^{-1}a$, we have 
$$(\sigma - q^{-1}\lambda)^{n}(ar) = q^{-n}a(\sigma-\lambda)^{n}(r)
$$ 
for all $r \in R$ and $n \ge 1$, and similarly for $ra$. Consequently, $aR_{\lambda} + R_{\lambda}a \subseteq R_{q^{-1}\lambda}$. In particular
\begin{equation}\label{eq:aT}
aT_i + T_ia \subseteq T_{i-1},
\quad\text{and hence}\quad
aT + Ta \subseteq T.
\end{equation}

\emph{Step 5: $T$ is semiprime.} Let $L \supseteq F$ be a finite field
extension containing all eigenvalues of $\sigma$, and let
$\Lambda \subseteq L^{*}$ be the finite set of all $\lambda$ with
$(R\otimes_F L)_{\lambda} \neq 0$. Then $R \otimes_F L = \bigoplus_{\lambda\in\Lambda}(R\otimes_F L)_{\lambda}$, and this decomposition is a grading by the subgroup $\Gamma$ of $L^{*}$
generated by $\Lambda$. Since $(R\otimes_F L)_{\lambda} = R_{\lambda}\otimes_F L$
for $\lambda \in F$, the algebra $T \otimes_F L$ is the identity component of
the induced grading by $\Gamma/H$, where $H = \Gamma \cap \{q^{m} : m\in \mathbb{Z}\}$. The
group $\Gamma/H$ is finitely generated and abelian, hence residually finite.
Therefore it has a subgroup $M$ of finite index such that the (finite)
image of $\Lambda \cup \{1\}$ in $\Gamma/H$ maps injectively into
$G = (\Gamma/H)/M$; thus $\lambda \in \Lambda$ has trivial image in $G$ if
and only if $\lambda \in H$. Hence $R \otimes_F L$ is graded by the finite
group $G$, with identity component $T \otimes_F L$.

Now let $x \in T$ satisfy $xTx = 0$. Since $x \otimes 1$ is homogeneous of
degree $1$, the left ideal $(R\otimes_F L)(x\otimes 1)$ is $G$-graded, and its
identity component $Tx \otimes_F L$ has square zero. By \cite{CR}, this left
ideal is nilpotent. Hence so is the left ideal $Rx$ of $R$, and the
semiprimeness of $R$ gives $Rx = 0$. Consequently, $x$ lies in the ideal
$\{y \in R : Ry = 0\}$, which has square zero, and therefore $x = 0$. This
shows that $T$ is semiprime.

\emph{Step 6: reduction to $T$.} The map $\sigma - 1$ is an algebraic
$1$-skew $\sigma$-derivation of both $R$ and $T$, and in both cases its
generalized $0$-eigenspace is $R_1$. By Theorem~\ref{thm:bg00}(1), applied
first to $R$ and then to $T$, both $R_1$ and $T$ have finite prime
dimension, and $R$ is PI whenever $R_1$ is.  As $R_1 \subseteq T$, it
suffices to show that $T$ is PI.

\emph{Step 7: the Jordan decomposition.} Setting $\sigma_s(r) = q^{i}r$ for
$r \in T_i$ defines an automorphism $\sigma_s$ of $T$ commuting with
$\sigma$, and $\sigma_u = \sigma\sigma_s^{-1}$ is a unipotent automorphism of
$T$, as $\sigma_u - 1 = q^{-i}(\sigma - q^{i})$ on $T_i$.  By \eqref{eq:aT},
$ar, ra \in T_{i-1}$ for $r \in T_i$, which gives, for all $r \in T$,
\begin{equation}\label{eq:jordan}
\begin{aligned}
\sigma_s(ar) &= q^{-1}a\sigma_s(r), & \sigma_s(ra) &= q^{-1}\sigma_s(r)a,\\
\sigma_u(ar) &= a\sigma_u(r), & \sigma_u(ra) &= \sigma_u(r)a.
\end{aligned}
\end{equation}

\emph{Step 8: the skew derivation $\de_u$.} For $r \in T$, set
$\de_u(r) = ar - \sigma_u(r)a$. By \eqref{eq:aT}, $\de_u(T) \subseteq T$, and
$\de_u$ is a $\sigma_u$-derivation of $T$. It follows from \eqref{eq:jordan}
that $\sigma_u\de_u = \de_u\sigma_u$ and $\de_u\sigma_s = q\,\sigma_s\de_u$.
Moreover, by \eqref{eq:jordan}, the operators $r \mapsto ar$ and
$r \mapsto \sigma_u(r)a$ commute, and the $j$-th power of the latter is
$r \mapsto \sigma_u^{j}(r)a^{j}$. Therefore
$$
\de_u^{\,m}(r) = \sum_{j=0}^{m}(-1)^{j}\binom{m}{j}a^{m-j}\sigma_u^{j}(r)a^{j},
$$
and since $a^{t} = 0$, each summand vanishes when $m \ge 2t-1$. Thus $\de_u$
is nilpotent.

\emph{Step 9: conclusion.} Since $\de_u\sigma_s = q\,\sigma_s\de_u$, we have
$\de_u(T_i) \subseteq T_{i-1}$ for all $i$, so $\de_u$ is homogeneous of
degree $-1$ with respect to the $\mathbb{Z}$-grading of $T$. Consequently, $T^{\de_u} = \bigoplus_{i\in\mathbb{Z}}(T^{\de_u} \cap T_i)$ is $\mathbb{Z}$-graded, with finite support and identity component
$T^{\de_u} \cap T_0 = (T^{\de_u})^{\sigma_s}$. If $r \in T^{\de_u}$ and
$\sigma_s(r) = r$, then $\sigma(r) = \sigma_u(r)$, and hence
$\de(r) = \de_u(r) = 0$. Thus $(T^{\de_u})^{\sigma_s} \subseteq R^{\de}$, and
therefore $(T^{\de_u})^{\sigma_s}$ is PI. Choose an integer $n$ with
$n > |i|$ for every $i$ in the support. Reducing degrees modulo $n$ yields a
grading of $T^{\de_u}$ by $\mathbb{Z}/n\mathbb{Z}$ with the same identity
component, so $T^{\de_u}$ is
PI by \cite[Corollary 9]{BC}. Summing up, $T$ is a semiprime algebra of finite prime
dimension, and $\de_u$ is a nilpotent $1$-skew $\sigma_u$-derivation of $T$,
with $\sigma_u$ algebraic, whose algebra of constants is PI. By
Theorem~\ref{thm:bg00}(3), $T$ is PI, and by Step~6 so is $R$.
\end{proof}

For a semiprime algebra $R$ let $Q = Q_s(R)$ be its symmetric Martindale ring
of quotients.  An automorphism $\sigma$ and a $\sigma$-derivation $\de$ of $R$
extend uniquely to $Q$, and we keep the notation.  Following Kharchenko \cite{Kh}, let $\Phi_{\sigma}$ be the set of $x \in Q$ such that $xr = \sigma(r)x$ for all $r \in R$.  The automorphism $\sigma$ is \emph{X-outer} if $\Phi_{\sigma} = 0$, and
\emph{X-inner} if $\Phi_{\sigma}$ contains a unit $b$ of $Q$, that is, if
$\sigma = \iota_b$ on $Q$, where $\iota_b(x) = bxb^{-1}$.  For semiprime $R$, $\sigma$ may be X-inner on some direct factors of $Q$ and X-outer on the others (Lemma \ref{lem:split}).

We work with $Z(R)^{\sigma}$, because $\de$ maps it into $\Phi_{\sigma}$
(Lemma \ref{lem:sigmanormal}).

\begin{lem}\label{lem:sigmanormal}
Let $\de$ be a $\sigma$-derivation of $R$ and let $z \in Z(R)^{\sigma}$.  Then
$a = \de(z)$ satisfies $ar = \sigma(r)a$ for every $r \in R$.  If $R$ is
semiprime, then $a \in \Phi_{\sigma}$, so that $a \neq 0$ forces $\sigma$ not
to be X-outer.
\end{lem}

\begin{proof}
Applying $\de$ to $zr = rz$ and using $\sigma(z) = z$ gives $\de(z)r + z\,\de(r) = \de(r)z + \sigma(r)\de(z)$, and $z$ is central, so $\de(z)r = \sigma(r)\de(z)$.  The second assertion
follows, as $R \subseteq Q$.
\end{proof}

\begin{lem}\label{lem:split}
Let $R$ be a semiprime PI algebra whose minimal prime ideals are
$P_1,\dots,P_k$, and let $\sigma$ be an automorphism and $\de$ a
$\sigma$-derivation of $R$, extended to $Q = Q_s(R)$.  Then
$Q =\prod_i Q_i$, where each $Q_i \cong Q_s(R/P_i)$ is a simple algebra,
finite dimensional over its centre $C_i$.  Let $e_i$ be the identity of
$Q_i$, so that $e_1,\dots,e_k$ are the primitive central idempotents of $Q$,
and let $e$ be the sum of those $e_i$ for which $\sigma(e_i) = e_i$ and
$\sigma|_{Q_i}$ is an inner automorphism of $Q_i$; for these $i$ we write
$e_i \le e$.  Then
\begin{enumerate}
\item $\sigma(e) = e$ and $\de(e) = 0$;
\item $\de(fx) = f\de(x)$ for all $x \in Q$ and every central idempotent
$f$ of $Q$ with $\sigma(f) = f$, in particular for $f = e$ and $f = e_i$
when $\sigma(e_i) = e_i$;
\item $\Phi_{\sigma} \subseteq eQ$, that is, $(1-e)\Phi_{\sigma} = 0$.
\end{enumerate}
\end{lem}

\begin{proof}
Let $J_i = \bigcap_{j \neq i}P_j$.
The sum $J_1 \oplus \dots \oplus J_k$ is an essential ideal of $R$, and $J_i$
is isomorphic to a nonzero ideal of the prime PI algebra $R/P_i$.  The
symmetric ring of quotients depends only on an essential ideal, so
$Q = \prod_i Q_s(R/P_i)$.  By Rowen's theorem \cite{Ro} every nonzero ideal
of $R/P_i$ meets its centre, so $Q_s(R/P_i)$ is the central closure of
$R/P_i$, which by Posner's theorem \cite{Po} is simple and finite dimensional
over its centre.  The $e_i$ are the primitive central idempotents of $Q$, so $\sigma$ permutes
them, say $\sigma(e_i) = e_{\pi(i)}$, and $\sigma(e) = e$.  If $f$ is a central
idempotent with $\sigma(f) = f$, then $\de(f) = \de(f^{2}) = \de(f)f + f\de(f)
= 2f\de(f)$.  Hence $f\de(f) = 0$ and $\de(f) = 0$, so
$\de(fx) = \de(f)x + \sigma(f)\de(x) = f\de(x)$.

Let $a \in \Phi_{\sigma}$.  For $x \in Q$ there is an essential ideal $D$ of
$R$ with $xD \subseteq R$, and then $axd = \sigma(xd)a = \sigma(x)ad$ for
$d \in D$, so $ax = \sigma(x)a$.  With $x = e_i$ this gives
$e_ia = e_{\pi(i)}a$, so $e_ia = e_ie_{\pi(i)}a = 0$ when $\pi(i) \neq i$.  If
$\pi(i) = i$, then $a_i = e_ia$ satisfies $a_iQ_i = \sigma(Q_i)a_i = Q_ia_i$,
a two-sided ideal of the simple ring $Q_i$.  If $a_i \neq 0$, this ideal is
all of $Q_i$.  Then $a_i$ is right invertible, hence invertible in the
artinian ring $Q_i$, and $\sigma|_{Q_i} = \iota_{a_i}$ is inner.  Thus $(1-e)a = 0$.
\end{proof}

\begin{lem}\label{lem:finiteorder}
Let $R$ be an algebra over a field $F$ of arbitrary characteristic, with
$Z(R)$ reduced, and let $\sigma$ be an algebraic automorphism of $R$ whose
minimal polynomial splits over $F$.  Then $\sigma$ has finite order on $Z(R)$.
\end{lem}

\begin{proof}
The centre $C = Z(R)$ is $\sigma$-stable, and
$C = \bigoplus_{\lambda}C_{\lambda}$ is the sum of its finitely many
generalized eigenspaces $C_{\lambda}$, with
$C_{\lambda} =\ker(\sigma-\lambda)^{k_{\lambda}}$.  Each nonzero $C_{\lambda}$ contains an eigenvector $c \neq 0$,
and then $\sigma(c^{\,n}) = \lambda^{n}c^{\,n}$ with $c^{\,n} \neq 0$, as $C$
is reduced.  As $\sigma$ is algebraic, the $\lambda^{n}$ are not pairwise
distinct, so $\lambda$ is a root of unity.  Let $n$ be a common multiple of
their orders and $\tau = \sigma^{n}$.  On $C_{\lambda}$ write
$\sigma = \lambda(1+N_{\lambda})$ with $N_{\lambda}^{k_{\lambda}} = 0$. Then
$\tau - 1 = (1+N_{\lambda})^{n} - 1$ is $N_{\lambda}$ times an operator
commuting with it, so $(\tau-1)^{k} = 0$ on $C$ for
$k = \max_{\lambda}k_{\lambda}$.  If $\operatorname{char}F = p > 0$ and
$p^{\,j} \ge k$, then $\tau^{\,p^{j}} - 1 = (\tau-1)^{p^{j}} = 0$ on $C$.  If
$\operatorname{char}F = 0$, suppose $\tau \neq 1$ on $C$.  Then
$\Delta = \tau - 1$ is a $\tau$-derivation of $C$ commuting with $\tau$.  Let
$l \ge 2$ be least with $\Delta^{l} = 0$ on $C$.  For $x \in C$ the Leibniz
formula for $\Delta^{2l-2}(x^{2}) = 0$ reduces to
$\binom{2l-2}{l-1}\tau^{l-1}\bigl(\Delta^{l-1}(x)\bigr)\Delta^{l-1}(x) = 0$.
As $\tau$ fixes $\Delta^{l-1}(x)$, we get $\Delta^{l-1}(x)^{2} = 0$, and
$\Delta^{l-1} = 0$ as $C$ is reduced, a contradiction.  So $\tau = 1$ on $C$.
\end{proof}

The case $\operatorname{char}F = 0$, $q = 1$ of the next theorem uses the
following fact about derivations of commutative algebras.

\begin{lem}\label{lem:topeigen}
Let $E$ be a reduced commutative algebra over a field $F$ of characteristic
$0$, $\partial$ a derivation of $E$, and $m \in F[t]$ a nonzero polynomial which
splits over $F$.  If $S \subseteq E$ is closed under multiplication and
$m(\partial)(S) = 0$, then $\partial(S) = 0$.
\end{lem}

\begin{proof}
For $\alpha \in F$ let $E_{\alpha} = \{x \in E : (\partial-\alpha)^{N}x = 0$ for some
$N\}$.  The sum of the $E_{\alpha}$ is direct, and the Leibniz rule
\[
(\partial-\alpha-\beta)^{N}(xy) = \sum_j\binom{N}{j}(\partial-\alpha)^{j}(x)\,
(\partial-\beta)^{N-j}(y)
\]
gives $E_{\alpha}E_{\beta} \subseteq E_{\alpha+\beta}$.
Let $\Lambda$ be the set of roots of $m$.  As $m$ splits,
$\ker m(\partial) \subseteq \bigoplus_{\alpha \in \Lambda}E_{\alpha}$.  Indeed,
write $m = \prod_{\alpha \in \Lambda}(t-\alpha)^{n_{\alpha}}$ and
$m_{\alpha} = m/(t-\alpha)^{n_{\alpha}}$.  The $m_{\alpha}$ have no common root,
so $1 = \sum_{\alpha}u_{\alpha}m_{\alpha}$ for some $u_{\alpha} \in F[t]$.  For
$x \in \ker m(\partial)$ this gives $x = \sum_{\alpha}u_{\alpha}(\partial)
m_{\alpha}(\partial)x$, and $(\partial-\alpha)^{n_{\alpha}}$ vanishes on the
$\alpha$-th term.

Let $0 \neq s \in S$ and write $s = \sum_{\alpha \in \Lambda}s_{\alpha}$ with
$s_{\alpha} \in E_{\alpha}$.  Let $\Lambda_s$ be the set of $\alpha$ with
$s_{\alpha} \neq 0$.  As $\mathbb{Q}$ is infinite, there is a
$\mathbb{Q}$-linear map $\ell \colon \mathbb{Q}\Lambda \to \mathbb{Q}$ on the
finite-dimensional $\mathbb{Q}$-span of $\Lambda$ which is injective on the
finite set $\Lambda \cup \{0\}$ (it must avoid finitely many hyperplanes).  Suppose that $\ell(\alpha) > 0$ for some $\alpha \in \Lambda_s$ and let
$\alpha^{*} \in \Lambda_s$ maximize $\ell$.  Then $s^{k}$ is the sum of the
products $s_{\beta_1}\cdots s_{\beta_k} \in E_{\beta_1+\dots+\beta_k}$, and
$\beta_1 + \dots + \beta_k = k\alpha^{*}$ forces
$\ell(\beta_j) = \ell(\alpha^{*})$, hence $\beta_j = \alpha^{*}$, for every
$j$.  So the component of $s^{k}$ in $E_{k\alpha^{*}}$ is
$s_{\alpha^{*}}^{\,k}$, which is nonzero as $E$ is reduced.  But
$s^{k} \in S$, so $s^{k} \in \bigoplus_{\alpha \in \Lambda}E_{\alpha}$, and
therefore $k\alpha^{*} \in \Lambda$ for all $k \ge 1$.  As
$\alpha^{*} \neq 0$ and $\operatorname{char}F = 0$, these are infinitely many
elements, a contradiction.  Applying this to $-\ell$ as well, we get
$\ell(\Lambda_s) = 0$, hence $\Lambda_s = \{0\}$ and $s \in E_0$.

Thus $0 \in \Lambda$ unless $S = \{0\}$, and, writing $m = t^{n}h$ with
$h(0) \neq 0$, we have $\partial^{n}(S) = 0$, since $s = s_0$ lies in
$\ker(\partial - 0)^{n}$ by the primary decomposition of $\ker m(\partial)$.  If $n \ge 2$,
then $x^{2} \in S$ for $x \in S$, so, as at the end of the proof of Lemma
\ref{lem:finiteorder} (with $\tau = \mathrm{id}$), $\partial^{2n-2}(x^{2}) = 0$
gives $\partial^{n-1}(S) = 0$.  By descending induction, $\partial(S) = 0$.
\end{proof}

\begin{thm}\label{thm:qskew}
Let $R$ be a semiprime algebra over a field $F$ with a $q$-skew
$\sigma$-derivation $\de$, where $\de$ and $\sigma$ are algebraic.  Assume that
$R^{\de}$ satisfies a polynomial identity and contains no infinite direct sum
of nonzero ideals.  Then $\GK(R) = \GK(R^{\de})$.
\end{thm}

\begin{proof}
First, $R$ has finite prime dimension.  The map $\de$
acts nilpotently on $R_0$ with $R_0^{\,\de} = R^{\de}$.  As $\sigma(R_0) = R_0$, as shown at the beginning of this section, the restriction of $\de$ is a
$q$-skew $\sigma$-derivation of $R_0$, and Lemma \ref{lem:primedim} shows
that $R_0$ has finite prime dimension.  So does $R$, by Theorem
\ref{thm:bg00}(1).

By Theorem \ref{thm:PIsemiprime} $R$ is PI, so it is semiprime PI of finite prime dimension
and, having finitely many minimal prime ideals, it satisfies $\GK(R) = \GK(Z)$, where $Z = Z(R)$, by Lemma \ref{lem:growthcentre}. Moreover, $Z$ is reduced.

\emph{Passing to $Z^{\sigma}$.}  We claim that $\sigma$ has finite order on
$Z$.  In characteristic $0$ put $\tau = \sigma$. In characteristic $p$ put
$\tau = \sigma^{\,p^{\ell}}$, whose eigenvalues are the $\lambda^{p^{\ell}}$
for the eigenvalues $\lambda$ of $\sigma$, separable over $F$ for $\ell$
large.  In both cases the minimal polynomial of $\tau$ splits over a finite
separable extension $K$ of $F$.  Since $Z(R \otimes_F K) = Z \otimes_F K$ is
reduced, Lemma \ref{lem:finiteorder} gives $\tau$, hence $\sigma$, finite
order on $Z$.  With $\nu$ the order of $\sigma$ on $Z$, every $z \in Z$ satisfies the
monic polynomial $\prod_{i<\nu}\bigl(X-\sigma^{i}(z)\bigr)$, whose coefficients
are symmetric in the orbit and hence $\sigma$-invariant.  Consequently, $Z$ is integral of bounded degree over $Z^{\sigma}$ and $\GK(Z) = \GK(Z^{\sigma})$.
So it suffices to show that $Z^{\sigma} \subseteq R^{\de}$, or, in
characteristic $p$, that $z^{p} \in R^{\de}$ for $z \in Z^{\sigma}$: then
$\GK(R) = \GK(Z^{\sigma}) \le \GK(R^{\de}) \le \GK(R)$, by Corollary
\ref{cor:integral} in the second case, as $Z^{\sigma}$ is integral of degree
at most $p$ over the subalgebra generated by the $z^{p}$.

\emph{The invariants of the centre.}  Let $Q = \prod_i Q_i$, the $e_i$ and $e$
be as in Lemma \ref{lem:split}.  For $z \in Z^{\sigma}$ the element $\de(z)$
lies in $\Phi_{\sigma}$ by Lemma \ref{lem:sigmanormal}, so $(1-e)\de(z) = 0$
by Lemma \ref{lem:split}(3), and
\begin{equation}\label{eq:star}\tag{$\ast$}
\de\bigl((1-e)z\bigr) \;=\; (1-e)\,\de(z) \;=\; 0 \qquad (z \in Z^{\sigma}).
\end{equation}
This settles the X-outer part.  It remains to treat the factors $Q_i$ with
$e_i \le e$.  Fix such an $i$ and a unit $b \in Q_i$ with $\sigma(x) = bxb^{-1}$
for $x \in Q_i$.  Then $\sigma(b) = b$, $\sigma$ fixes $C_i$ pointwise, and
$\de(Q_i) \subseteq Q_i$.  The map $\de' = b^{-1}\de$ is an ordinary derivation
of $Q_i$, since
$\de'(xy) = b^{-1}\de(x)y + b^{-1}(bxb^{-1})\de(y) = \de'(x)y + x\de'(y)$, and
in particular $\de'(C_i) \subseteq C_i$.  Put $w = b^{-1}\de(b)$.  Expanding
$\de(bxb^{-1}) = q\,b\,\de(x)\,b^{-1}$ with the help of
$\de(b^{-1}) = -b^{-1}\de(b)b^{-1}$ and cancelling $b$ gives
\begin{equation}\label{eq:starstar}\tag{$\ast\ast$}
(1-q)\,\de(x) \;=\; \sigma(x)w - wx \qquad (x \in Q_i).
\end{equation}

\emph{The case $q \neq 1$.}  By \eqref{eq:starstar}, $\de$ is an inner
$\sigma$-derivation of $Q_i$, so it vanishes on the central elements fixed by
$\sigma$, that is on $C_i$.  Hence $\de(e_iz) = 0$ for $z \in Z$, and with
\eqref{eq:star}, $\de(z) = \sum_{e_i \le e}\de(e_iz) + \de\bigl((1-e)z\bigr) = 0$ for
$z \in Z^{\sigma}$.  Thus $Z^{\sigma} \subseteq R^{\de}$.

\emph{The case $\operatorname{char}F = p > 0$.}  For $z \in Z$ and $e_i \le e$
the element $e_iz$ is central in $Q_i$, so
$\de'\bigl((e_iz)^{p}\bigr) = p(e_iz)^{p-1}\de'(e_iz) = 0$ and
$\de\bigl((e_iz)^{p}\bigr) = 0$.  For $z \in Z^{\sigma}$ also
$(1-e)z \in Q^{\de}$ by \eqref{eq:star}, and $Q^{\de}$ is a subring, because
$\de(xy) = \de(x)y + \sigma(x)\de(y)$.  The $e_i$ and $1-e$ are orthogonal
central idempotents, so
\[
z^{p} = \sum_{e_i \le e}(e_iz)^{p} + \bigl((1-e)z\bigr)^{p}
\in Q^{\de} \cap R = R^{\de}.
\]

\emph{The case $\operatorname{char}F = 0$ and $q = 1$.}  Let $K$ be a finite
extension of $F$ over which the minimal polynomial $m$ of $\de$ splits, put
$R_K = R \otimes_F K$ and extend $\sigma$ and $\de$ $K$-linearly.  In
characteristic $0$ the algebra $R_K$ is again semiprime, as $K/F$ is
separable, and PI.  It has finite
prime dimension, since every ideal of $R_K$ is an $R$-sub-bimodule of
$R_K \cong R^{\,[K:F]}$ and $R$ has finite uniform dimension as an
$R \otimes_F R^{op}$-module.  So Lemmas \ref{lem:sigmanormal} and
\ref{lem:split}, and hence \eqref{eq:star}, apply to $R_K$.  As
$Z^{\sigma} \subseteq Z(R_K)^{\sigma}$ and $R_K^{\,\de} \cap R = R^{\de}$, it
is enough to show $Z(R_K)^{\sigma} \subseteq R_K^{\,\de}$.  We may therefore
assume $K = F$, that is, that $m$ splits over $F$.

Fix $e_i \le e$ and $b$ as above.  For $q = 1$, \eqref{eq:starstar} reads
$\sigma(x)w = wx$, that is $x\,b^{-1}w = b^{-1}w\,x$ for $x \in Q_i$, so
$\beta = b^{-1}w$ lies in $C_i$ and $\de(b) = bw = b^{2}\beta$.  Let
$E = C_i[b]$, a commutative subalgebra of $Q_i$.  The automorphism
$\sigma = \iota_b$ is the identity on $E$, so for $x,y \in E$ we have
$\de(xy) = \de(x)y + x\de(y)$.  Moreover $\de(C_i) = b\de'(C_i) \subseteq E$
and $\de(b) \in E$, so $\de(E) \subseteq E$, and $\partial = \de|_{E}$ is a derivation
of the commutative $F$-algebra $E$.

Let $S = e_iZ$.  It is closed under multiplication and lies in
$C_i \subseteq E$, and since $\de(e_ix) = e_i\de(x)$ by Lemma \ref{lem:split}(2),
$m(\partial)(e_iz) = e_i\,m(\de)(z) = 0$ for $z \in Z$.
As $E$ is finite dimensional over the field $C_i$, it is artinian, hence
noetherian, and as $\operatorname{char}F = 0$, the nilradical $\mathcal{N}$
of $E$ is $\partial$-stable \cite[Proposition 14.2.3]{MR}.  Thus $\partial$ induces a derivation $\bar\partial$ of the reduced algebra
$\bar E = E/\mathcal{N}$, and the image $\bar S$ of $S$ satisfies
$m(\bar\partial)(\bar S) = 0$.  Lemma \ref{lem:topeigen} gives $\bar\partial(\bar S) = 0$,
that is, $\de(e_iz) \in \mathcal{N}$ for every $z \in Z$.  Now
$\de(e_iz) = b\,\de'(e_iz)$ with $\de'(e_iz) \in C_i$.  As $Q_i$ is finite
dimensional over $C_i$, the unit $b$ is a root of a polynomial over $C_i$ with
nonzero constant term, so $b^{-1} \in E$.  Hence $\de'(e_iz) = b^{-1}\de(e_iz)$ is a nilpotent element of the field $C_i$.  Therefore $\de'(e_iz) = 0$ and $\de(e_iz) = 0$.

Together with \eqref{eq:star} this gives
$\de(z) = \sum_{e_i \le e}\de(e_iz) + \de\bigl((1-e)z\bigr) = 0$ for every
$z \in Z^{\sigma}$, that is, $Z^{\sigma} \subseteq R^{\de}$.
\end{proof}

\begin{cor}\label{cor:semiprime}
Let $R$ be a semiprime algebra over a field $F$ with an algebraic derivation $d$.  Assume that
$R^{d}$ satisfies a polynomial identity and contains no infinite direct sum
of nonzero ideals.  Then $\GK(R) = \GK(R^{d})$.
\end{cor}

\begin{proof}
Take in Theorem \ref{thm:qskew} $\sigma = \mathrm{id}$ and $q = 1$. 
\end{proof}

\section{Examples}\label{sec:examples}

In this section we show that the hypotheses of the Main Theorem cannot be
substantially weakened.

Without semiprimeness of $R$ and finiteness of the prime dimension of
$R^{d}$ the conclusion of the Main Theorem fails, already for derivations,
and the difference of the two GK dimensions can be made arbitrarily large while $d^{2} = 0$.  The following construction, given here with several variables, makes this
precise.

\begin{example}\label{ex:biggap}
Fix $r \ge 1$, let $P = F[x_1,\dots,x_r]$ and let $W = P \ast F[y]$ be the
free product of $P$ and $F[y]$.  It has a basis consisting of $1$ and the
\emph{alternating words} $u_1u_2\cdots u_l$ $(l \ge 1)$, whose factors are
alternately nonconstant monomials of $P$ and positive powers of $y$, such as
$x_1x_2^{3}\,y\,x_2$.  Let
\[
A_r = \left\{ \begin{pmatrix} f & g \\ 0 & f \end{pmatrix} :
f \in P,\ g \in W \right\},
\]
and write $(f,g)$ for such a matrix, so that
$(f,g)(f',g') = (ff',\, fg' + gf')$ and $1 = (1,0)$.  Let $b = (0,y)$ and
$d = \ad_{b}$.  Then $d^{2} = 0$, and $d(f,g) = (0,\, yf - fy)$ vanishes exactly when
$f \in F$, since no nonconstant element of $P$ commutes with $y$ in $W$.  Hence
\[
A_r^{\,d} = F \oplus N, \qquad N = \{(0,g) : g \in W\}, \qquad N^{2} = 0,
\]
so every affine subalgebra of $A_r^{\,d}$ is finite dimensional and
$\GK(A_r^{\,d}) = 0$.  The algebra $A_r^{\,d}$ is commutative, hence PI, but
its prime dimension is infinite: since $(\lambda,g)(0,h) = (0,\lambda h)$,
every subspace of $N$ is an ideal of $A_r^{\,d}$.

On the other hand $\GK(A_r) = 2r$.  For the lower bound let $U$ be the span of $1$,
$(x_1,0),\dots,(x_r,0)$ and $b$.  For monomials $u,v$ of $P$ we have
$(u,0)\,b\,(v,0) = (0,uyv)$, and the words
$uyv$ are linearly independent for distinct pairs $(u,v)$, being distinct
alternating words of the free product.  The number of pairs with
$\deg u + \deg v \le n-1$ is $\binom{n-1+2r}{2r}$, the number of monomials of
degree at most $n-1$ in $2r$ variables, so
$\dim_F U^{n} \ge \binom{n-1+2r}{2r}$, which is of order $n^{2r}$.  For the
upper bound, $A_r$ is not finitely generated, so consider any finite
dimensional subspace $U_0$.  It lies in
$\{(u,0) : \deg u \le c\} + \{(0,g) : g \in S\}$ for some $c$ and
some finite-dimensional $S \subseteq W$, where $u$ runs over polynomials.  A
product with two factors from $N$ vanishes, so $U_0^{\,n}$ is spanned by the
$(u,0)$ with $u$ a monomial of degree at most $cn$ and the $(0, ugv)$ with $g$
in a basis of $S$ and $u,v$ monomials with $\deg u + \deg v \le cn$.  Hence
$\dim_F U_0^{\,n} \le \binom{cn+r}{r} + (\dim S)\binom{cn+2r}{2r}$, which
grows like $n^{2r}$, and therefore $\GK(A_r) \le 2r$.

Thus the difference $\GK(A_r) - \GK(A_r^{\,d}) = 2r$ can be arbitrarily large, although $d^{2} = 0$.  With infinitely many
variables $x_i$ it is infinite, as for the resulting algebra $A_{\infty}$
the constants $A_{\infty}^{\,d} = F \oplus N$ still have GK dimension $0$.  The hypotheses that fail here are the semiprimeness of $A_r$, since $N$ is a nonzero nilpotent
ideal of $A_r$, and the finiteness of the prime dimension of $A_r^{\,d}$,
as noted above.  Thus this example shows only that these two hypotheses cannot
be dropped together.\qed
\end{example}

Semiprimeness of $R$ together with finite prime dimension of $R^{d}$ does
not suffice either: for a primitive algebra
$R \supseteq A_r$ and the inner derivation $d = \ad_b$, with $b$ as in
Example~\ref{ex:biggap}, the gap $\GK(R) - \GK(R^{d})$ can be at least as
large.  This follows from a general construction.  Let $A$ be an algebra
acting faithfully on a space with a fixed basis, under the mild conditions of
Lemma \ref{lem:finitary}.  We adjoin to $A$ the finitary matrices
with respect to this basis.  The resulting algebra $R$ is primitive, and
$\GK(R) \le 2\GK(A)$.  Moreover, for $b \in A$ and $d = \ad_b$ we have
$\GK(R^{d}) \le 2\GK(C_A(b))$.

\begin{lem}\label{lem:finitary}
Let $A$ be a subalgebra of $\End_F(V)$ containing $\mathrm{id}_V$, where $V$
is an infinite-dimensional vector space with basis $\mathcal{M}$.  For $\alpha,\beta \in \mathcal{M}$ let
$E_{\alpha\beta}$ be the matrix unit sending $\beta$ to $\alpha$ and the other
basis elements to $0$, and let $\mathcal{E}$ be their span.  Suppose that
\begin{enumerate}[label=(A\arabic*)]
\item for each $a \in A$ and $\beta \in \mathcal{M}$, the $\beta$-coordinate of
$ax$ is nonzero for only finitely many $x \in \mathcal{M}$;
\item no nonzero element of $A$ has finite rank.
\end{enumerate}
Then $R = A + \mathcal{E}$ is a primitive algebra which is not PI, $\mathcal{E}$
is an ideal of $R$, $R = A \oplus \mathcal{E}$, and
$\GK(A) \le \GK(R) \le 2\GK(A)$.  Let $b \in A$ and $d = \ad_{b}$.  Then
$R^{d} = C_A(b) \oplus \mathcal{E}^{d}$, where $C_A(b)$ is the centralizer of
$b$ in $A$, and
$$
\GK\bigl(C_A(b)\bigr) \;\le\; \GK\bigl(R^{d}\bigr) \;\le\;
2\,\GK\bigl(C_A(b)\bigr) .
$$
Also $\GK(R^{d}) = 0$ if $C_A(b)$ is locally finite.

Condition (A1) holds when $V = A$, with $A$ acting by left multiplication,
$A$ is graded by an abelian group with finite-dimensional homogeneous
components, and $\mathcal{M}$ consists of homogeneous elements.
\end{lem}

\begin{proof}
We begin by fixing the notation.  For $\beta \in \mathcal{M}$ let $c_{\beta}\colon V \to F$ be
the $\beta$-coordinate functional, so that $x = \sum_{\beta}c_{\beta}(x)\beta$
for $x \in V$, with only finitely many nonzero terms.  Let $V'$ be the span of
the $c_{\beta}$, that is, the space of finite linear combinations of coordinate
functionals.  For $v \in V$ and a linear functional $\psi$ on $V$ let
$v \otimes \psi \in \End_F(V)$ be the operator $x \mapsto \psi(x)v$.  It is bilinear in $(v,\psi)$, and
$\alpha \otimes c_{\beta} = E_{\alpha\beta}$ for $\alpha,\beta \in \mathcal{M}$.
So for $v = \sum_{\alpha}\lambda_{\alpha}\alpha$ and
$\psi = \sum_{\beta}\mu_{\beta}c_{\beta}$, both sums finite, we get
$v \otimes \psi = \sum_{\alpha,\beta}\lambda_{\alpha}\mu_{\beta}E_{\alpha\beta}$.
Hence
\begin{equation}\label{eq:spanE}
\mathcal{E} \;=\; \operatorname{span}\{E_{\alpha\beta}\}
\;=\; \operatorname{span}\{\,v \otimes \psi : v \in V,\ \psi \in V'\,\} .
\end{equation}

\emph{Step 1: $\mathcal{E}$ is an ideal of $R$.}  Let $a \in \End_F(V)$,
$v,v' \in V$ and $\psi,\psi' \in V'$.  For every $x \in V$,
$$
\bigl(a(v \otimes \psi)\bigr)(x) = a\bigl(\psi(x)v\bigr) = \psi(x)\,av ,
\qquad
\bigl((v \otimes \psi)a\bigr)(x) = \psi(ax)\,v ,
$$
so that
\begin{equation}\label{eq:rules}
\begin{gathered}
a(v \otimes \psi) = av \otimes \psi, \qquad
(v \otimes \psi)a = v \otimes (\psi \circ a),\\
(v \otimes \psi)(v' \otimes \psi') = \psi(v')\,v \otimes \psi' .
\end{gathered}
\end{equation}
These identities hold for every $a \in \End_F(V)$.  Now let $a \in A$.  The
first and the third operators lie in $\mathcal{E}$ by \eqref{eq:spanE}.  For the second
we need $\psi \circ a \in V'$, and by linearity it suffices to take
$\psi = c_{\beta}$.  We claim that $c_{\beta} \circ a = \sum_{x \in \mathcal{M}} c_{\beta}(ax)\,c_x$.  The sum is finite by (A1), and both sides take the value $c_{\beta}(ax')$ at
every $x' \in \mathcal{M}$, so they agree.  Thus $A\mathcal{E}$,
$\mathcal{E}A$ and $\mathcal{E}\mathcal{E}$ lie in $\mathcal{E}$.  Therefore $R$ is a subalgebra of $\End_F(V)$, and $\mathcal{E}$ is an ideal of $R$.

\emph{Step 2: $R$ is primitive and not PI.}  As $R \subseteq \End_F(V)$, the
space $V$ is a faithful $R$-module.  It is irreducible.  Indeed, let
$0 \neq x \in V$ and $v \in V$.  Choose $\beta$ with $c_{\beta}(x) \neq 0$.
Then $(v \otimes c_{\beta})(x) = c_{\beta}(x)v$, so $v \in Rx$.  Hence $R$ is
primitive.  For every $n$ choose distinct $\alpha_1,\dots,\alpha_n \in
\mathcal{M}$, which is possible as $\mathcal{M}$ is infinite.  The matrix units
$E_{\alpha_i\alpha_j}$ span a subalgebra of $\mathcal{E}$ isomorphic to
$M_n(F)$, which satisfies no polynomial identity of degree less than $2n$.  Hence $R$ is not PI.

\emph{Step 3: Decomposition  $R = A \oplus \mathcal{E}$.}  Every element of $\mathcal{E}$ has
finite rank, so $A \cap \mathcal{E} = 0$ by (A2).

\emph{Step 4: Computation of constants.}  Let $b \in A$ and $d = \ad_{b}$.  Then
$[b,A] \subseteq A$, and $[b,\mathcal{E}] \subseteq \mathcal{E}$ as
$\mathcal{E}$ is an ideal.  Let $a \in A$ and $\varphi \in \mathcal{E}$ with
$a + \varphi \in R^{d}$.  Then $[b,a] = -[b,\varphi]$ lies in
$A \cap \mathcal{E} = 0$, so $a \in C_A(b)$ and $\varphi \in \mathcal{E}^{d}$.
Hence $R^{d} = C_A(b) \oplus \mathcal{E}^{d}$.

\emph{Step 5: Computation of the growth.}  The lower bounds hold because $A \subseteq R$ and
$C_A(b) \subseteq R^{d}$.  For the upper bounds let either $B = A$ and
$\mathcal{E}' = \mathcal{E}$, or $B = C_A(b)$ and
$\mathcal{E}' = \mathcal{E}^{d}$, so that $B \oplus \mathcal{E}'$ is $R$ or
$R^{d}$.  Let $S$ be an affine subalgebra of $B \oplus \mathcal{E}'$.  Its
generators lie in $Y + \operatorname{span}\{\varphi_1,\dots,\varphi_s\}$ for
some finite-dimensional subspace $Y \ni 1$ of $B$ and some
$\varphi_1,\dots,\varphi_s \in \mathcal{E}'$.  By \eqref{eq:spanE} write
$\varphi_i = \sum_l x_{il} \otimes \psi_{il}$, a finite sum with
$x_{il} \in V$ and $\psi_{il} \in V'$, and put
$X = \operatorname{span}\{x_{il}\}$ and $X' = \operatorname{span}\{\psi_{il}\}$.
Let $S_n$ be the span of the products of at most $n$ factors from
$Y \cup \{\varphi_1,\dots,\varphi_s\}$.  A product without a factor
$\varphi_i$ lies in $Y^{n}$.  A product with some factor $\varphi_i$ has the
form $u_0\varphi_{i_1}u_1\varphi_{i_2}\cdots\varphi_{i_k}u_k$ with $k \ge 1$ and
$u_0,\dots,u_k \in Y^{n}$.  Expanding each $\varphi_i$ and using the rules
\eqref{eq:rules}, it is a sum of terms $\lambda\,(u_0x) \otimes (\psi \circ u_k)$ with
$\lambda \in F$, $x \in X$ and $\psi \in X'$.  For instance, let $k = 2$ and
take one summand $x \otimes \psi$ of $\varphi_{i_1}$ and one summand
$x' \otimes \chi$ of $\varphi_{i_2}$.  By \eqref{eq:rules},
$$
\begin{aligned}
u_0(x \otimes \psi)u_1(x' \otimes \chi)u_2
&= \bigl(u_0x \otimes (\psi \circ u_1)\bigr)\bigl(x' \otimes (\chi \circ u_2)\bigr)\\
&= \psi(u_1x')\,(u_0x) \otimes (\chi \circ u_2) .
\end{aligned}
$$
In general the scalar is
$\lambda = \psi_1(u_1x_2)\,\psi_2(u_2x_3)\cdots\psi_{k-1}(u_{k-1}x_k)$, where
$x_j \otimes \psi_j$ is the chosen summand of $\varphi_{i_j}$ (so $\lambda = 1$ if $k = 1$), and then
$x = x_1$ and $\psi = \psi_k$.  Therefore $S_n \subseteq Y^{n} + (Y^{n}X) \otimes (X'Y^{n})$, where $Y^{n}X$ is spanned by the $ux$ and $X'Y^{n}$ by the $\psi \circ u$, for
$u \in Y^{n}$, $x \in X$ and $\psi \in X'$.  Hence $\dim S_n \le \dim Y^{n} + (\dim X)(\dim X')(\dim Y^{n})^{2}$.  For every $\varepsilon > 0$ we have $\dim Y^{n} \le n^{\GK(B)+\varepsilon}$ for
all large $n$, so $\GK(S) \le 2\GK(B)$.  As $S$ was arbitrary,
$\GK(B \oplus \mathcal{E}') \le 2\GK(B)$.  If $B$ is locally finite, then $Y$
lies in a finite-dimensional subalgebra $B_0$ of $B$, so $Y^{n} \subseteq B_0$
for all $n$, and $S$ is finite dimensional.
\medskip

Finally, we consider the graded case. Let $V = A$ with left multiplication,
$A = \bigoplus_{\gamma \in \Gamma}A_{\gamma}$ with $\dim A_{\gamma} < \infty$,
and $\mathcal{M}$ homogeneous.  By linearity it suffices to check (A1) for
$a \in A_{\eta}$.  For $x \in \mathcal{M}$ of degree $\kappa$ we have $ax \in A_{\kappa+\eta}$.  Hence, if $\beta$ has degree $\mu$, then $c_{\beta}(ax) \neq 0$ forces $\kappa = \mu - \eta$, and only finitely
many elements of $\mathcal{M}$ have degree $\mu - \eta$.
\end{proof}

\begin{example}\label{ex:primitive}
Let $A_r$ and $b$ be as in Example \ref{ex:biggap}.  Grade $A_r$ by
multidegree in the $x_i$ and $y$, with values in the free abelian group
$\Gamma$ on these variables: $(A_r)_{\gamma}$ consists of the $(f,g)$ with
$f \in P_{\gamma}$ and $g \in W_{\gamma}$.  The relations $x_ix_j = x_jx_i$ are homogeneous, and in
$(f,g)(f',g') = (ff',\,fg' + gf')$ all three products have degree
$\gamma + \gamma'$, so this is a grading.  Each multidegree is realized by
finitely many words, so the components are finite dimensional, also with
infinitely many $x_i$.  Let $\mathcal{M}$ be the homogeneous basis of the
elements $(m,0)$ and $(0,w)$, where $m$ is a monomial of $P$ and $w$ is $1$
or an alternating word of $W$, and let $A_r$ act on $V = A_r$ by left
multiplication.  Then (A1) holds by the last assertion of Lemma
\ref{lem:finitary}.  For (A2), a nonzero $(f,g)$ sends $(f',0)$ to
$(ff', gf')$.  If $f \neq 0$, then $f' \mapsto ff'$ is injective on $P$; if
$f = 0$, then $g \neq 0$ and $f' \mapsto gf'$ is injective, as the free
product $W$ of domains is a domain.  So $(f,g)$ has infinite rank.

Let $R = A_r + \mathcal{E}$ and $d = \ad_{b}$.  As $b^{2} = 0$, we have
$d^{\,3} = 0$.  By Example \ref{ex:biggap}, $C_{A_r}(b) = A_r^{\,d} = F \oplus N$ with
$N^{2} = 0$, which is locally finite.  Lemma \ref{lem:finitary} gives that $R$
is primitive and not PI, and $\GK(R^{d}) = 0 < 2r \le \GK(A_r) \le \GK(R)$.  With infinitely many variables $x_i$ one gets $\GK(R) = \infty$.

We show that a nonzero subspace $J_P$ of $R^{d}$ lies in every nonzero ideal
of $R^{d}$, so any two nonzero ideals of $R^{d}$ intersect nontrivially.  We
use the notation $c_{\beta}$, $V'$ and $v \otimes \psi$ of the proof of Lemma
\ref{lem:finitary}.  Here $b(f,g) = (0,yf)$, so
$\ker b = \{(0,g) : g \in W\}$ and
$\operatorname{im} b = \{(0,yf) : f \in P\}$, which is spanned by the basis
elements $(0,ym)$, $m$ a monomial.  Let $\Psi$ be the span of the $c_{\beta}$
with $\beta$ not of this form, and $\Psi_P \subseteq \Psi$ the span of the
$c_{(m,0)}$.  Then $\psi \in V'$ vanishes on $\operatorname{im} b$ if and only
if $\psi \in \Psi$.  Every $\rho \in R^{d}$ commutes with $b$, so
$\rho(\ker b) \subseteq \ker b$.  We record five facts.

(a) If $v \in \ker b$ and $\psi \in \Psi$, then $v \otimes \psi \in R^{d}$,
since $b(v \otimes \psi) = bv \otimes \psi = 0$ and
$(v \otimes \psi)b = v \otimes (\psi \circ b) = 0$.

(b) Let $c_m$ be the coordinate functionals of $P$ on the monomial basis.
For $h \in P$ and a monomial $m'$ let $h \otimes c_{m'}$ be the operator
$f \mapsto c_{m'}(f)h$ of $P$; finite sums of such operators are called
finitary matrices on $P$.  For $g \in W$ let $g_y \in P$ be the polynomial with
$c_m(g_y) = c_{(0,ym)}(0,g)$ for every monomial $m$, so that
$(yf)_y = f$ for $f \in P$.
For a finitary matrix $X$ let $\varphi_X(f,g) = (Xf,\,yX(g_y))$.  If
$X = h \otimes c_{m'}$, then
$\varphi_X = (h,0) \otimes c_{(m',0)} + (0,yh) \otimes c_{(0,ym')}$, as both
sides send $(f,g)$ to $\bigl(c_{m'}(f)h,\ c_{m'}(g_y)\,yh\bigr)$.  By
linearity and \eqref{eq:spanE}, $\varphi_X \in \mathcal{E}$ for every $X$.
Since $(yf)_y = f$,
$$
b\varphi_X(f,g) = (0,\,yXf) = \varphi_X(0,\,yf) = \varphi_X b(f,g),
$$
so $\varphi_X \in \mathcal{E}^{d} \subseteq R^{d}$, and $\varphi_X$ acts on
$\operatorname{im} b$ by $(0,yf) \mapsto (0,yXf)$.  Now let $\chi \in \Psi$
and write $\chi = \psi + \chi'$, where
$\psi = \sum_m \mu_m c_{(m,0)} \in \Psi_P$ and $\chi'$ is a combination of the
$c_{(0,w)}$ with $w$ not of the form $ym$.  As $\varphi_X$ takes values in
the span of the $(m,0)$ and the $(0,ym)$, we get $\chi' \circ \varphi_X = 0$
and $(\chi \circ \varphi_X)(f,g) = \tilde\psi(Xf)$, where
$\tilde\psi = \sum_m \mu_m c_m$.  In particular, for $X = h \otimes c_{m'}$
and $v \in V$, \eqref{eq:rules} gives
$$
(v \otimes \chi)\varphi_X = v \otimes (\chi \circ \varphi_X)
= \tilde\psi(h)\, v \otimes c_{(m',0)} .
$$

(c) Let $J$ be the span of the $v \otimes \psi$ with $v \in \ker b$ and
$\psi \in \Psi$.  If $t \in \mathcal{E}^{d}$ has image in $\ker b$, then
$t \in J$.  Indeed, write $t = \sum_j u_j \otimes \chi_j$ with the $u_j$ and
the $\chi_j$ linearly independent.  The $u_j$ span the image of $t$, so
$u_j \in \ker b$ and $bt = 0$.  Hence
$0 = tb = \sum_j u_j \otimes (\chi_j \circ b)$, which gives $\chi_j \in \Psi$.

(d) If $u \in \ker b \setminus \operatorname{im} b$, there is $\psi \in \Psi$
with $\psi(u) \neq 0$, since $u$ has a nonzero coordinate $c_{(0,w)}(u)$ with
$w$ not of the form $ym$.

(e) If $t = \sum_k v_k \otimes \psi_k$ and $t = \sum_j u_j \otimes \chi_j$
with the $u_j$ and the $\chi_j$ linearly independent, then the $u_j$ span
$t(V) \subseteq \operatorname{span}\{v_k\}$ and the $\chi_j$ span
$\{\lambda \circ t\} \subseteq \operatorname{span}\{\psi_k\}$, $\lambda$
ranging over linear functionals.  So every element of $J$, $J_0$ or $J_P$
below has such a representation with factors in the corresponding spaces.

Put
\[
J_0 = \operatorname{span}\{v \otimes \psi : v \in \operatorname{im} b,\
\psi \in \Psi\}, \qquad J_P = \operatorname{span}\{v \otimes \psi :
v \in \operatorname{im} b,\ \psi \in \Psi_P\}.
\]
Both lie in $R^{d}$ by (a),
as $\operatorname{im} b \subseteq \ker b$, and $J_P \neq 0$.  Let $I$ be a
nonzero ideal of $R^{d}$ and $0 \neq \rho \in I$.

\emph{Step 1: $I \cap J \neq 0$.}  If $\rho(u) \neq 0$ for some
$u \in \ker b$, then $\rho(u \otimes c_{(1,0)}) = \rho(u) \otimes c_{(1,0)}$
is a nonzero element of $I \cap J$.  Otherwise $\rho(\ker b) = 0$, so
$b(\rho x) = \rho(bx) = 0$ for every $x$, that is
$\rho V \subseteq \ker b$, and there is $f \in P$ with $\rho(f,0) \neq 0$.
For $X = f \otimes c_{1}$ we have $\varphi_X(1,0) = (f,0)$, so
$t = \rho\varphi_X$ is a nonzero element of $I \cap \mathcal{E}^{d}$ with
image in $\ker b$, and $t \in J$ by (c).

\emph{Step 2: $I \cap J_0 \neq 0$.}  Let $0 \neq t = \sum_j u_j \otimes \chi_j
\in I \cap J$ with the $u_j \in \ker b$ and the $\chi_j \in \Psi$ linearly
independent.  If all $u_j$ lie in $\operatorname{im} b$, then $t \in J_0$.
Otherwise choose, by (d), $\psi \in \Psi$ with $\psi(u_j) \neq 0$ for some
$j$, and $0 \neq v \in \operatorname{im} b$.  Then
$(v \otimes \psi)t = v \otimes \sum_j \psi(u_j)\chi_j$ is a nonzero element
of $I \cap J_0$.

\emph{Step 3: $J_P \subseteq I$.}  Let $0 \neq t = \sum_j v_j \otimes \chi_j
\in I \cap J_0$ with $v_j = (0,yf_j)$, the $f_j$ and the $\chi_j$ linearly
independent.  For $g \in P$ there is a finitary $X$ with $Xf_1 = g$ and
$Xf_j = 0$ for $j \ge 2$, and then $\varphi_X t = (0,yg) \otimes \chi_1$.  So
$v \otimes \chi_1 \in I$ for every $v \in \operatorname{im} b$.  If the
$\Psi_P$-component $\psi$ of $\chi_1$ is nonzero, choose $h \in P$ with
$\tilde\psi(h) = 1$; by (b), $(v \otimes \chi_1)\varphi_{h \otimes c_{m'}} =
v \otimes c_{(m',0)} \in I$ for every monomial $m'$, so $J_P \subseteq I$.
Otherwise $\chi_1(0,w_0) \neq 0$ for some $w_0$ not of the form $ym$.  As
$(0,w_0) \in \ker b$, (a) gives $(0,w_0) \otimes c_{(1,0)} \in R^{d}$, and
$(v \otimes \chi_1)\bigl((0,w_0) \otimes c_{(1,0)}\bigr)
= \chi_1(0,w_0)\,v \otimes c_{(1,0)}$.  So $v \otimes c_{(1,0)} \in I$ for
every $v \in \operatorname{im} b$, and the first case applies to $c_{(1,0)}$
in place of $\chi_1$.

Hence every nonzero ideal of $R^{d}$ contains $J_P \neq 0$.  So $R$ and
$R^{d}$ satisfy every hypothesis of the Main Theorem except the polynomial
identity on $R^{d}$.\qed
\end{example}

\section{Questions}\label{sec:questions}

By Example~\ref{ex:biggap}, the assumptions of Corollary~\ref{cor:semiprime}
that $R$ is semiprime and that $R^{d}$ has finite prime dimension cannot be
dropped together.  Can either of them be dropped separately? More precisely, let $d$ be an algebraic derivation of an algebra $R$ such that $R^{d}$
satisfies a polynomial identity.  Is $\GK(R) = \GK(R^{d})$ in either of the following cases:

\begin{enumerate}[label=(\alph*)]
\item $R$ is semiprime and $R^{d}$ contains an infinite
direct sum of nonzero ideals.
\item $R$ is not necessarily semiprime, and $R^{d}$ contains no infinite
direct sum of nonzero ideals ?
\end{enumerate}

In case (a), $R$ satisfies a polynomial identity by \cite{CL}.

\end{document}